%
%
%
%
%
%
\documentstyle{amsppt}
\magnification \magstep1
\parskip 11pt
\parindent .3in 
\pagewidth{5.2in} 
\pageheight{7.2 in}

%
%
\def \bl{\vskip 11pt} 

\def \ni{\noindent}
\def \P{\bold{P}}
\def \Q{\bold{Q}}

\def \C{\bold{C}} 
\def \O{\Cal{O}}
\def \lra{\longrightarrow}
%

\nologo
\vskip 15pt
\centerline {\bf A Remark on Projective Embeddings of
Varieties with }
\centerline{\bf  Non-Negative Cotangent Bundles}
\vskip 10pt
\centerline{\smc Lawrence EIN \footnote{Partially
supported by NSF Grant DMS 96-22540}}
\centerline{\smc Bo ILIC  \footnote{Partially supported
by an NSERC postdoctoral fellowship}}
\centerline{\smc Robert LAZARSFELD  \footnote{Partially
supported by NSF 
 Grant DMS 97-13149}}
\vskip 15pt
\centerline{\it Dedicated to the memory of Michael
Schneider}

\vskip 15pt
\ni{\bf Introduction.}

The purpose of this note is to establish an elementary
but somewhat unexpected bound on the degrees of
projective embeddings of varieties with numerically
effective cotangent bundles. 

In recent years, there has been interest in
understanding the geometry of  complex projective
varieties whose tangent or cotangent bundles satisfy
various positivity properties.  In this note, we shall
be concerned with smooth complex projective varieties
$X$ satisfying the following non-negativity property:
\

\ni (NCB).\ \ \ \    {\it The cotangent bundle 
$\Omega^1_X$ of $X$ is numerically effective
$($nef$)$.}

\ni By definition, the condition means that the Serre
line bundle
$\O_{\P(\Omega^1_X)}(1)$ on the projectivization
$\P(\Omega^1_X)$ is numerically effective, or
equivalently that for any non-constant map $ \nu : C
\lra X$ from a smooth curve $C$ to $X$, any quotient
bundle of
$\nu^* \Omega^1_X$ has non-negative degree. Property
(NCB) is satisfied, for example, by smooth subvarieties
of abelian varieties, by varieties uniformized by the
ball or other irreducible Hermitian symmetric spaces
(cf. \cite{Mok}, \S1), and by products and submanifolds
thereof. 

Our result is that if $X$ satisfies (NCB), then  the
degree of $X$ in any projective embedding must grow
essentially exponentially in the dimension of $X$.
Specifically, given a positive integer $n$, define 
$$\delta(n) = 2 ^{ [  \sqrt{n}   ] },$$
 where as usual $[x]$ denotes the integer part of $x$. 
\proclaim {\bf Theorem}  Let  $X$ be a smooth
projective variety of dimension $n$ which satisfies
Property (NCB), and let 
$$f : X \lra \P^n$$ be any finite surjective mapping.
Then $\deg(f) \ge
\delta(n)$. In particular, the degree of
$X$ in any projective embedding $X \subset \P^r$ must
be at least
$\delta(n)$. \endproclaim

\ni We suspect that these statements are not optimal,
and that there should be genuinely exponential, or even
factorial, bounds on the degree. It would be
interesting to prove results along these lines. In
another direction, it seems natural to wonder  whether
similar degree bounds hold also for varieties whose
universal covers are e.g. bounded Stein domains.   More
philosophically, these results suggest that the
complexity of the projective geometry associated to
varieties satisfying (NCB) grows exponentially with
their dimension. It would be interesting to know if one
could make this viewpoint precise, and whether it has
any other manifestations.

The proof of the Theorem requires only a few lines, and
in fact the two ingredients that enter into the
argument are at least implicitly quite well known.  One
simply notes that the hypothesis (NCB) forces the
presence of points where the derivative of
$f$ drops rank substantially, and that this in turn
leads to a lower bound on $\deg(f)$. Nonetheless, the
conclusion came as something of a surprise to us: while
linear bounds on the degree are very familiar (eg.
\cite{GL}, Theorem 2), the existence of essentially
exponential statements seems to have been overlooked. 

The third author had the opportunity to discuss some of
these matters with Michael Schneider  about  a year
before his
 death, and as always Michael was enthusiastic and
encouraging. We hope therefore that the present note
might not be out of place in this volume dedicated to
his memory. Schneider contributed a lot to algebraic
geometry on both a personal and a professional level,
and he will be greatly missed. 

The proof of the main result occupies in \S 1. Some
applications and variants appear in \S 2. We are
grateful to D. Burns and N. Mok for some valuable
discussions. 

\vskip 15pt

\ni {\bf  \S 1. Proof  of the Theorem.}

We start with a lemma on degrees and singularities of
branched coverings. It was suggested by some examples
of Flenner and Ran alluded to  in
\cite{Ran 2}. 
\proclaim{\bf Lemma 1.1}  Let $f: X \lra Y$ be a finite
surjective map of smooth complex varieties of dimension
$n$. Fix a point $x \in X$, let $y = f(x) \in Y$, and
denote by
$e_f(x)$ the local degree of $f$ at $x$, i.e. the
multiplicity of $x$ in its fibre $f^{-1} f(x)$. 
Suppose that derivative $df_x : T_xX \lra T_yY$ of $f$
at $x$ has rank $n - k$. Then
$e_f(x) \ge 2^k$, and consequently $\deg(f) \ge 2^k$.
\endproclaim

\ni {\smc Proof}. By hypothesis, the co-derivative
$df_x^* : T^*_yY \lra T^*_xX$ has a $k$-dimensional
kernel. Denoting by $m_x \subset \O_xX$ and
$m_y \subset \O_yY$ the maximal ideals of $x$ and $y$
respectively,  we can therefore choose a system of
parameters
$u_1, \dots , u_n \in m_y$ in such a way that $f^* u_1,
\dots , f^* u_k \in m_x^2$. Now
$$e_f(x) = \dim_{\C} \O_xX / f^* m_y , $$ i.e. $e_f(x)$
is alternatively the intersection multiplicity at $x$
of the (germs of) divisors defined by the $f^* u_i$. On
the other hand,  it is well known (cf
\cite{Fult, 12.4}) that this intersection multiplicity
is at least the product of the multiplicities
$\text{ord}_x f^* u_i$ of the individual divisors.
Since by construction 
$\text{ord}_x f^* u_i \ge 2$ for $1 \le i \le k$, the
stated lower bound on
$e_f(x)$ follows. The inequality on $\deg(f)$ is then a
consequence the fact that for fixed $y \in Y$, 
$$\sum_{f(x) = y}  e_f(x) = \deg(f).  \qed $$

The plan is to apply the Lemma to branched coverings of
projective space. The following well-known fact, which
we include for the convenience of the reader,  will let
us apply theorems on degeneracy loci to guarantee the
existence of singularities.
\proclaim{\bf Lemma 1.2}  Let $X$ be a projective
variety, and let $E$ and
$F$ be vector bundles on $X$.  If $E$ is nef and $F$ is
ample, then $E
\otimes F$ is ample.
\endproclaim
\ni {\smc Sketch of  Proof.} The statement is a
consequence of Kleiman's criterion (cf.
\cite{Hart}) that the nef cone is the closure of the
ample cone, and the argument is most easily stated
using the language of vector bundles twisted by
$\Q$-divisors, as in \cite{Myka}.  First, one verifies
the statement when
$F$ is a line bundle, or more generally an ample
$Q$-divisor: we leave this to the reader. Next, fix an
ample line bundle
$H$ on $X$.  Since $E$ is nef, it follows that
$E(\frac{1}{N}H)$ is ample for any $N > 0$, and since
$F$ is ample, $F(-\frac{1}{N}H)$ is ample for $N
\gg 0$. Therefore
$E \otimes F = E(\frac{1}{N}H) \otimes  F(-
\frac{1}{N}H)$ is ample. \qed

Now we turn to the 

\ni {\smc Proof of the Theorem}. Assume  that $X$ is
smooth projective variety of dimension  $n$ whose
cotangent bundle $\Omega_X^1$ is nef, and suppose given
a branched covering $f : X \lra \P^n$. Let
$$S_i(f) = \left \{ x \in X \mid \text{rank} \  df_x
\le n - i \right \}.$$ This is an algebraic subset of
$X$ whose expected dimension is $n - i^2$ (cf.
\cite{Fult}, Chapter 14). In particular, setting $k = [
\sqrt{n} ]$,
$S_k(f)$ has non-negative postulated dimension. The
asserted bound on
$\deg(f)$ will follow from Lemma 1 as soon as we show
that $S_k(f) \ne
\emptyset$. But this is a consequence of \cite{FL1} or
\cite{L, \S2} or
\cite{FL2}. In fact, since the tangent bundle $T\P^n$
(and hence also $f^* T \P^n$) is ample, the hypothesis
(NCB) implies by Lemma 2 that 
$\Omega^1_X \otimes f^* T\P^n$ is an ample vector
bundle on $X$. The cited results then guarantee that
the vector bundle map $df : TX \lra f^* T\P^n$ must
actually drop rank whenever it is dimensionally
predicted to do so. Finally, given an embedding $X
\subset \P$ of $X$ into some projective space, we get
by projection  a branched covering $f : X \lra \P^n$
whose degree is the degree of $X$ in $\P$, and so
$\deg(X) \ge \delta(n)$.
\qed

\ni{\bf Remark}. Given a smooth variety $X$ with nef
cotangent bundle, and an ample line bundle $L$ on $X$
which is generated by its global sections, the theorem
is equivalent to the assertion that $\int c_1(L)^n \ge
\delta(n)$. It is perhaps worth noting that this bound
can fail if $L$ is not globally generated. For example,
fixing $n$, let $C$ be a  smooth curve of genus $g \gg
n$ which carries no $g^1_n$, and let $X = Sym^n(C)$ be
the
$n^{\text{th}}$ symmetric product of $C$. The
Abel-Jacobi map $X \lra Jac^n(C)$ is an embedding, so
$X$ satisfies (NCB). On the other hand, upon choosing a
base-point $P \in C$, 
 $Sym^{n-1}(C)$ embeds as a divisor in  $X$ (via $D
\mapsto D + P$), and the corresponding line bundle $L =
\O_X(Sym^{n-1}(C))$ is ample (cf. \cite{FL1,
\S2}). But
$\int_X c_1(L)^n = 1$, as one sees from the fact that
there is a unique effective divisor of degree $n$
containing $n$ given points of $C$.

\vskip 15 pt
\ni{\bf \S 2. Applications and Variants.}

\vskip 10pt

We begin with a simple application of the Theorem:
\def \G {\bold G}

\proclaim{Corollary 2.1} Let $A$ be an abelian variety
of dimension $m$, and let $X \subset A$ be a smooth
subvariety of dimension $n$. Assume that
$X$ is of general type. Then the top self-intersection
of the canonical bundle of $X$ satisfies the inequality:
$$\int c_1(\O_X(K_X))^n \ge \delta(n).$$ \endproclaim

\ni{\smc Proof of Corollary.} 	The embedding $X\subset
A$ gives rise to a Gauss mapping
$\gamma : X \lra \G$ of
$X$ into the Grassmannian
$\G = \G(n,m)$ of $n$-dimensional subspaces of $T_0A$,
which is generically finite since $X$ is of general
type (cf. \cite{Mori, \S 3}).  A theorem of Ran
\cite{Ran1} implies that then $\gamma$ is actually
finite. On the other hand, the Pl\"ucker line bundle
$\O_\G(1)$ on $\G$ pulls back to the canonical bundle
on $X$. Therefore the canonical bundle
$\O_X(K_X)$ is ample and globally generated. But $X$ --
like any submanifold of $A$ --
 satisfies Property (NCB), and  the desired inequality
then follows from the Theorem.
\qed

We next prove a variant of the Theorem for certain
smooth subvarieties of projective space:
\proclaim{Proposition 2.2} Let $X \subset  \P^{n+e} =
\P$ be a smooth subvariety of projective space having
dimension $n$ and codimension $e$, and denote by $N =
N_{X/\P}$ the normal bundle to $X$ in $\P$.  If $N(-1)$
is ample, then 
$$\deg(X) \ge \min \left \{ 2^e, \delta(n) \right \}.$$
 \endproclaim
\ni Recall that at least when $X$ spans $\P$, 
the hypothesis on $N(-1)$ is equivalent
to requiring that every hyperplane tangent to $X$ be
tangent at only finitely many points. Note that we do
not assume here that
$X$ satisfies (NCB). Observe also that if $e^2 \le n$,
then the stated  bound
$\deg(X) \ge 2^e$ is best possible for a complete
intersection of quadrics.

\ni {\smc Proof of Proposition 2.2.} Fix a linear space
$L^{e-1}$ disjoint from $X$, and project from $L$ to
get a finite mapping $f : X \lra \P^n$. Setting $k =
\min \{e , [\sqrt{n}] \}$, we will show that the
singularity locus $S_k(f)$ appearing in the proof of
the Theorem is non-empty, and then the result will
follow as above from Lemma 1.1. To this end, recalling
that  $\P^{n+e} - L$ is the total space of
$\O_{\P^n}(1)^{\oplus e}$, one finds the exact sequence
of bundles on $X$:
$$ 0 \lra \O_X(1)^{\oplus e} \lra T \P^{n+e} | X \lra
f^* T \P^n \lra 0.$$ Combining this with the sequence
$$0 \lra TX \lra T \P^{n+e} | X \lra N_{X/\P} \lra 0,$$
we arrive at a mapping of vector bundles
$$u : \O_X(1)^{\oplus e} \lra N_{X/\P}$$ on $X$ whose
degeneracy loci are the same as the degeneracy loci of
the derivative
$df : TX \lra f^* T \P^n$. Since the bundle
$N_{X/\P}(-1)$ is ample, the results cited in the proof
of Theorem 1 imply that $S_{k}(u) \ne
\emptyset$, as desired. \qed

\ni{\bf Exercise 2.3.} Suppose that $X \subset \P^{n+e}
$ is a smooth subvariety having the property that for
some $x \in X$ the embedded tangent space $T_xX \subset
\P^{n+e}$ meets $X$ at only finitely many points (so
that in particular $e \ge n$). Then
$deg(X) \ge 2^n$.

An argument  similar to the one proving Proposition 2.2
also leads to the following generalization of the Main
Theorem:
\proclaim{\bf Proposition 2.4} Let $X$ be a smooth
variety of dimension
$n$ satisfying (NCB), and let
$E$ be an ample vector bundle of rank $e$ on $X$ which
is generated by its global sections. Then
$$\int_X  s_n(E) \ge \min \left \{ 2^n , \delta(n+e -1)
\right \},$$ where
$s_n(E)$ denotes the $n^{\text{th}}$ Segre class of
$E$.\endproclaim
\demo{Outline of Proof} In brief, consider the
projective bundle
$\pi : \P(E) \lra X$, and fix a general subspace
$V \subset H^0(X, E)$ of dimension $n + e$ generating
$E$. This gives rise to a finite mapping $f : \P(E)
\lra \P(V) = \P^{n + e - 1}$ whose degree is equal to
$\int s_n(E)$. Setting $k = \min \{n,  [\sqrt{n+ e -
1}] \}$, it is enough as above to show that the
singularity locus $S_k(f)$ is non-empty. To this end,
let $M$ be the vector bundle of rank $n$ on $X$ defined
by the exact sequence
$$0 \lra M \lra V \otimes_\C \O_X \lra E \lra 0, \tag *
$$ the homomorphism on the right being the canonical
evaluation map. Now $f$ factors through the embedding
$\P(E) \subset \P(V \otimes_\C \O_X) = X \times \P(V) 
$ determined  by (*), and as in the proof of the
Proposition, the degeneracy loci of $df$ coincide with
those of the resulting vector bundle map $$u :
\pi ^* TX \lra N_{\P(E) / X \times \P(V)}.$$ But
$\P(E)$ is cut out in $X \times \P(V)$ by a section of
$pr_1^*M^* \otimes pr_2^*
\O_{\P(V)}(1)$, and consequently $N_{\P(E) / X \times
\P(V)} = \pi^*M^*
\otimes
\O_{\P(E)}(1)$, which by Lemma 2 is ample thanks to the
amplitude of $E$ and the fact that
$M^*$ is globally generated. As $\pi^* \Omega^1_X$ is
nef by assumption, it follows that 
$\pi^* \Omega^1_X \otimes  (\pi^*M^* \otimes
\O_{\P(E)}(1))$ is ample. But then
\cite{FL1} or the other references cited above
guarantee that $u$ must actually drop rank whenever it
is dimensionally predicted to do so. \qed
\enddemo

\ni {\bf Remark 2.5}. The inequalities established in
this note all spring via Lemma 1.1 from producing
singularities of a branched covering of projective
space. It would be interesting to know whether one can
recover or improve these statements by applying
positivity theorems to some well-chosen Chern class
calculations. It is natural to wonder in particular
whether the inequalities of
\cite{BSS} might not be relevant here.

\vskip 10pt

\ni 

\vskip 20pt
\bl
\settabs\+University of Illinois at Chicago and now is
the time  \cr 
\+ Lawrence EIN \cr
\+ Department of Mathematics \cr
\+ University of Illinois at Chicago \cr
\+ 851 South Morgan St., M/C  249 \cr
\+ Chicago, IL  60607-7045 \cr

\vskip 7pt

\+ Bo ILIC  \cr
\+ Department of Mathematics \cr
\+ U.C.L.A. \cr
\+ Los Angeles, CA  90095 \cr

\vskip 7 pt

\+ Robert LAZARSFELD \cr
\+ Department of Mathematics \cr
\+ University of Michigan \cr
\+ Ann Arbor, MI  48109 \cr

\end